\documentclass[11pt]{article}
\usepackage{latexsym,amsfonts}
\usepackage{amssymb,amsmath}
\usepackage{enumitem}
\usepackage{graphics}
\usepackage[demo]{graphicx}
\usepackage{epsfig}
\usepackage{pstricks}
\usepackage[normal]{subfigure}
\usepackage[latin1]{inputenc}
\usepackage[english,francais]{babel}
\usepackage{relsize,exscale}
\usepackage{makeidx}
\usepackage{enumitem}
\usepackage{amsfonts,amssymb,amsmath}
\usepackage{graphicx}
\usepackage{color}
\usepackage{multirow}
\usepackage{mathrsfs}
\usepackage[normalem]{ulem}
\usepackage{cancel}
\usepackage{bbm}
\usepackage{empheq}
\newenvironment{prooff}{{\it Proof :}}{\hfill\rule{2mm}{2mm}\vskip3mm\par}
\newtheorem{theorem}{Theorem}[section]
\newtheorem{lemma}[theorem]{Lemma}

\newtheorem{e-definition}[theorem]{Definition\rm}
\newtheorem{remark}{\it Remark\/}

\usepackage{color}
\definecolor{dred}{rgb}{0.92,0,0}
\definecolor{dgreen}{rgb}{0,0.92,0}
\definecolor{dblue}{rgb}{0,0,0.92}
\definecolor{dyellow}{rgb}{0.95,0.95,0}

\newcommand{\Z}{\mathbb{Z}}
\newcommand{\R}{\mathbb{R}}
\newcommand{\N}{\mathbb{N}}
\def\D{\displaystyle}

\newcommand{\hs}{\hspace{0.1cm}}

\newcommand{\sa}{\\ [0.2cm]}
\usepackage{multirow}
\DeclareMathOperator*{\esssup}{ess\,sup}
\graphicspath{
{./Figures/}
{Figures/}
{./}
}
\title{Enhancing Interpolation and Approximation Error Estimates\\
 Using a Novel Taylor-like Formula}
\author{Jo\"el Chaskalovic \thanks{D'Alembert,
Sorbonne University, Paris, France, (Email:  jch1826@gmail.com)}
\qquad
Franck Assous
\thanks{
Ariel University, Ariel, Israel, (Email: \emph{(corresp.)} assous@ariel.ac.il)}
\qquad
}
\date{}
\begin{document}
\maketitle
\selectlanguage{english}
\begin{abstract}
\noindent  In this paper, we present an approach to enhance interpolation and approximation error estimates.
Based on a previously derived first-order Taylor-like formula, we demonstrate its applicability in improving the $P_1$-interpolation error estimate. Following the same 
principles, we also develop a novel numerical scheme for the heat equation that yields a better error estimate compared to the classical implicit finite differences scheme.
\end{abstract}

\noindent {\em keywords}: Taylor's theorem, Taylor-like formula, interpolation error, Finite differences scheme, Error estimate, Heat equation.
\section{Introduction}\label{intro}

\noindent Even today, improving the accuracy of approximation continues to pose a challenge in the field of numerical analysis. In this context, we recently derived a novel 
first-order Taylor-like formula in \cite{arXiv_First_Order}. The aim of this formula was to obtain a reduced remainder compared to the one obtained with the traditional 
Taylor formula. This result was achieved by ``transferring'' the numerical weight of the  remainder to the main part of the expansion. Extension to second order were also proposed in \cite{ChAs2024}, in the same spirit as what can be found in \cite{DrACxx} for quadrature formula \sa
\noindent From a mathematical perspective, the root of these issues can be found in Rolle's theorem and in Lagrange and Taylor's theorems, (see for instance 
\cite{Atki88}, \cite{BuFa11}), basically due to an unknown point in the remainder of the Taylor expansion. As a result, most of the error estimates concentrate on the 
{\em asymptotic behavior} of the error, for instance, in finite element methods, as the mesh size tends to zero.\sa
\noindent In this context, various approaches have been proposed to identify ways of enhancing the approximation accuracy. For instance, in the realm of numerical integration, 
we refer the reader to \cite{Barnett_Dragomir}, \cite{Cerone} or \cite{Dragomir_Sofo}, along with the references cited therein. From a different perspective, because of a the 
lack of information, heuristic methods have been explored, basically based on a probabilistic approach, see for instance \cite{Abdulle}, \cite{AsCh2014}, \cite{Hennig}, 
\cite{Oates} or \cite{ArXiv_JCH}, \cite{CMAM2} and \cite{ChAs20}. This makes possible to compare different numerical approaches, for instance for a fixed mesh size, see 
\cite{MMA2021}.\sa
\noindent However, accurately estimating the upper bounds of error estimates and developing methods to enhance these upper bounds remain significant challenges. In this context, 
we proposed in \cite{ChAs2023} a refined first-order expansion formula in $\R^n$,  to obtain a reduced remainder compared to the one obtained by usual Taylor's formula, and we 
investigated some related properties. Similar problems have been considered in the past years, and are often referred to as the perturbed (or corrected) quadrature rules, see 
for instance \cite{Cerone} or \cite{Dragomir_Sofo}.  In other instances,  authors obtained in \cite{Chen01}, \cite{DrWa97} or \cite{MaPU00} the trapezoid inequality by the 
difference between the supremum and the infimum bounds of the first derivative.\sa
\noindent In this paper, we address error estimates related to interpolation and approximation based on our new Taylor-like formula. We will demonstrate that the interpolation 
error estimate can be enhanced, even within the broader class of Sobolev spaces $W^{1,1}$. Additionally, we explore how finite differences schemes can be 
improved (in terms of precision), by using our approach.\sa
\noindent The paper is organized as follows. In Section \ref{A-new-first-order}, we recall the first-order Taylor-like formula proposed in \cite{arXiv_First_Order}. Then, we 
consider classical applications of the new Taylor expansion. In Section \ref{Improving-P1-interpolation}, we show how it can improve the $P_1$-interpolation error estimate. 
In Section \ref{Section-newsecond}, we use it to derive a new numerical scheme for the heat equation. Concluding remarks follow.
\section{A new first order expansion formula}\label{A-new-first-order}
\noindent To begin with, let us recall the first order Taylor-like formula we derived in \cite{arXiv_First_Order}. To this end, we consider an integer $n$ in $\N^*$, $(a,b) \in \mathbb{R}^{2}$, $a < b$, and a function $f \in \mathcal{C}^{2}(]a,b[)$. We have
\begin{theorem}
Let $f$ be a real mapping defined on $[a,b]$ which belongs to $\mathcal{C}^{2}([a,b])$, such that: $\forall x \in [a,b], -\infty < m_2 \leqslant f''(x) \leqslant M_2 < +\infty$. \sa
Then, we have the following first order expansion:
\begin{equation}\label{T1}
f(b) = f(a) + (b-a)\left(\frac{f'(a) + f'(b)}{2n} + \frac{1}{n}\sum \limits_{k=1}^{n-1} f'\left(a + k\frac{(b-a)}{n}\right)\right) + (b-a)\epsilon_{a,n+1}^{(1)}(b),
\end{equation}
where :
\begin{equation}\label{Remainder_TL}
\D|\epsilon_{a,n+1}^{(1)}(b)| \leqslant \frac{(b-a)}{8n}(M_2-m_2).
\end{equation}
\end{theorem}
This result is optimal in the sense that the weights involved in the linear combination of $f'$ at the equally spaced points $\D a + k\frac{(b-a)}{n}$ guarantee 
the remainder $\epsilon_{a,n+1}^{(1)}(b)$ to be minimum.\sa
The particular case of (\ref{T1}) when $n=1$ will be also considered in the sequel. In this case, we have:
\begin{equation}\label{T2}
\D f(b) = f(a) + (b-a)\left(\frac{f'(a)+f'(b)}{2}\right) + (b-a)\epsilon_{a,2}(b)\,.
\end{equation}
Above, the remainder $\epsilon_{a,2}(b)$ satisfies the following inequality:
$$
\frac{(b-a)}{8}(m_2-M_2)\leqslant \epsilon_{a,2}(b) \leqslant \frac{(b-a)}{8}(M_2-m_2).
$$
This result can be compared with the well-known first order Taylor's formula \cite{Taylor}, that we recall here for completeness. With the same notations, we have:
$$
f(b) = f(a) + (b-a)f'(a) + (b-a)\epsilon_{a,1}(b),
$$
where:
$$
\D \frac{(b-a)}{2}m_2 \leqslant \epsilon_{a,1}(b) \leqslant \frac{(b-a)}{2}M_2.
$$

\noindent In the following of the paper, we will consider two applications where formulas (\ref{T1}) and (\ref{T2}) will be involved. In Section \ref{Improving-P1-interpolation} we will show how (\ref{T1}) can be used to improve the $P_1$-interpolation error estimate. In Section \ref{Heat}, we will apply (\ref{T2}) to develop a new numerical scheme for the heat equation, with a smaller upper bound in the approximation error estimate 
compared to the one associated with the classical implicit finite differences scheme.
\section{Improving the $W^{1,1}$ interpolation error estimate}\label{Improving-P1-interpolation}
\noindent The first application we have in mind is related to the interpolation error estimate. Now, numerous partial differential equations are not well posed for any 
integer $m$ in $H^m(\Omega)$ but in a more general class of Sobolev spaces, namely, $W^{m,p}(\Omega), (m,p)\in\N^{*2}$, and in particular in the space $W^{1,1}(\Omega)$, where $\Omega$ denotes a given non empty open domain in $\R^d$.\sa
It is for example the case of the Laplace equation with a given right-hand side $f \in L^{p}(\Omega), (p\ne 2)$. Indeed, in that case, the solution $u$ to the associated 
variational formulation belongs to $W^{1,p}(\Omega), (p\ne 2)$, if the domain $\Omega$ is regular enough. Other cases may be found for example  in \cite{BeSc72}, 
\cite{BrNi83} and \cite{AtEP01}. \sa
In the sequel, we consider for $\Omega = ]0,1[$, the functional framework based on the Sobolev space $W^{1,1}(]0,1[)$ defined by:
$$
\D W^{1,1}(]0,1[) = \left\{u:\,]0,1[\rightarrow\R, \frac{}{}u \in L^1(]0,1[)\,;\, u' \in L^1(]0,1[)\right\},
$$
where $u'$ is the weak derivative of $u$ in $L^1(]0,1[)$, see \cite{Brezis}.\sa\
\noindent Using this framework, we will derive in this section a new interpolation error estimate based on the new Taylor-like formula (\ref{T1}). More precisely, 
we consider a given real function $u$ defined on the interval $[0,1]$ which belongs to $C^2([0,1])\subset W^{1,1}(]0,1[)$. \sa
We also introduce a mesh on $[0,1]$ defined by: $(x_i)_{i=0,\dots,N+1}$ such that $x_0=0$ and $x_{N+1}=1$. Moreover, we define the mesh size $h$ by: 
$\D h=\max_{i=0,\dots,N}h_i$, where $h_i=x_{i+1}-x_i, (i=0,\dots,N)$. \sa
Finally, we consider the $P_1$-interpolation polynomial $u_{I}$ of $u$ which satisfies: 
\begin{eqnarray}
\forall i \in \{0,\dots,N+1\}, \, u_{I}(x_i)=u(x_i), & & \nonumber\\[0.2cm]
\forall x \in [x_i,x_{i+1}], u_{I} \in P_1([x_i,x_{i+1}]), & & \nonumber
\end{eqnarray}
where $P_1([x_i,x_{i+1}])$ is the space of polynomials of degree at most 1 defined on $[x_i,x_{i+1}]$.
We also introduce the  following notations: For any $u\in W^{1,1}(]0,1[)$, we denote by $\|.\|_{1,1}$ the standard norm defined by:
$$
\D \|u\|_{1,1} = \|u\|_{0,1} + \|u'\|_{0,1},
$$
where the norm $\|u\|_{0,1}$ is defined by:
$$
\D \|u\|_{0,1} = \int_{0}^{1}|u(x)|\,dx.
$$
We first derive the classical interpolation error estimate based on the standard Taylor formula:
\begin{lemma}
Let u be in $C^2([0,1])$ and $u_I$ the corresponding $P_1$-interpolation polynomial. Then, the standard Taylor formula leads to the following interpolation error estimate:
\begin{equation}\label{Error_Estim_T}
\D \|u-u_I\|_{1,1} \leq (h+h^2)\|u''\|_{\infty},
\end{equation}
where $\D\|u''\|_{\infty}=\esssup_{x\in[0, 1]}|u''(x)|.
$
\end{lemma}
\begin{prooff}
We recall the classical first order Taylor formula \cite{arXiv_First_Order} given by:
\begin{equation}\label{Taylor}
u(x_{i+1}) = u(x_i) + h_i u'(x_i) + h_i \epsilon^{(T)},
\end{equation}
with:
\begin{equation}\label{Reste_Taylor}
\D|\epsilon^{(T)}| \leq \frac{h_i}{2}\|u''\|_{\infty}.
\end{equation}
- We begin by evaluating the $L^{1}$-norm of the derivative, that is $\|u'-u'_{I}\|_{0,1}$. \sa
We have:
$$
\D \|u'-u'_{I}\|_{0,1} = \int_{0}^{1}|u'(x)-u'_{I}(x)|\,dx = \sum_{i=0}^{N}\int_{x_i}^{x_{i+1}}|u'(x)-u'_{I}(x)|\,dx.
$$
Then, given that $u'_I$ is constant on $[x_i,x_{i+1}]$,  by the help of (\ref{Taylor}), we get:
$$
\D \forall \in [x_i,x_{i+1}]: u'_{I}(x) = \frac{u(x_{i+1})-u(x_i)}{h_i} = u'(x_i) + \epsilon^{(T)}.
$$
As a consequence, using (\ref{Reste_Taylor}) and Fubini's theorem \cite{Brezis}, we can derive the following sequence of inequalities:
\begin{eqnarray}
\D \int_{x_i}^{x_{i+1}}|u'(x)-u'_{I}(x)|\,dx & = & \int_{x_i}^{x_{i+1}}|u'(x)-u'(x_i) - \epsilon^{(T)}|\,dx \label{IN1}\\[0.2cm]
\D & \leq & \int_{x_i}^{x_{i+1}}|u'(x)-u'(x_i)|\, dx +\frac{h_i^2}{2}\|u''\|_{\infty} \nonumber \\[0.2cm]
\D & \leq & \int_{x_i}^{x_{i+1}}\bigg|\int_{x_i}^{x}u''(t)dt\bigg|\,dx +\frac{h_i^2}{2}\|u''\|_{\infty} \nonumber \\[0.2cm]
\D & \leq & \int_{x_i}^{x_{i+1}}\int_{x_i}^{x}|u''(t)|dt\, dx +\frac{h_i^2}{2}\|u''\|_{\infty}  \nonumber \\[0.2cm]
\D & \leq & \int_{x_i}^{x_{i+1}}|u''(t)|\left(\int_{t}^{x_{i+1}}\!\!dx\right)dt +\frac{h_i^2}{2} \|u''\|_{\infty}  \nonumber \\[0.2cm]
\D & \leq & \int_{x_i}^{x_{i+1}}(x_{i+1}-t)|u''(t)|\,dt +\frac{h_i^2}{2}\|u''\|_{\infty} = h_i^2 \|u''\|_{\infty}. \label{IN6}
\end{eqnarray}
By summing over $i$ between 0 and $N$, and using that $\D\sum_{i=0}^{N}h_i=1$, we finally get:
\begin{equation}\label{N1_u1}
\D \|u'-u'_{I}\|_{0,1} \leq \bigg(\sum_{i=0}^{N}h_i^2\bigg)\|u''\|_{\infty} \leq h\|u''\|_{\infty}\,.
\end{equation}
\sa
\noindent - Let us now evaluate the $L^{1}$-norm $\|u-u_{I}\|_{0,1}$. First of all, we remark that we have, for all $x \in  [x_i,x_{i+1}]$,
\begin{equation}\label{Ineq0}
 |u(x)-u_{I}(x)| = \bigg| \int_{x_i}^{x}\big(u'(t)-u_I'(t)\big)\,dt \bigg|.
\end{equation}
Then,
\begin{equation}\label{Ineq00}
\D |u(x)-u_{I}(x)| \leq \int_{x_i}^{x}|u'(t)-u'_I(t)|\,dt \leq \int_{x_i}^{x_{i+1}}|u'(t)-u'_I(t)|\,dt.
\end{equation}
So, by using inequalities (\ref{IN1})-(\ref{IN6}), we get that
$$
\D \forall \in [x_i,x_{i+1}],  \, |u(x)-u_{I}(x)| \leq h_i^2 \|u''\|_{\infty}.
$$
It remains now to integrate this inequality on $[x_i,x_{i+1}]$ to obtain that
$$
\D \int_{x_i}^{x_{i+1}}|u(x)-u_{I}(x)| \leq h_i^3 \|u''\|_{\infty},
$$
and summing over all values of $i$ between $0$ and $N$, this yields
\begin{equation}\label{N1_u2}
\D \|u-u_{I}\|_{0,1} \leq h^2\|u''\|_{\infty}.
\end{equation}
Finally, by combining inequalities (\ref{N1_u1}) and (\ref{N1_u2}), we get the interpolation error estimate (\ref{Error_Estim_T}).
\end{prooff}
Now, let us derive the interpolation error estimate obtained by using the Taylor-like formula (\ref{T1}), instead of the standard Taylor formula.
\begin{theorem}\label{Thmfirst}
Let u be in $C^2([0,1])$ and $u_I$ the corresponding $P_1$-interpolation polynomial. Then, the Taylor-like formula (\ref{T1}) leads to the following interpolation error estimate:
\begin{equation}\label{Error_Estim_TL}
\D \|u-u_I\|_{1,1} \leq \frac{(h+h^2)}{2}\|u''\|_{\infty} + \frac{(h+h^2)}{8n}(M_2-m_2),  \forall n\in\N^{*}.
\end{equation}
In particular, when $n$ goes to $+\infty$, we get:
\begin{equation}\label{Error_Estim_TL_Asymp}
\D \|u-u_I\|_{1,1} \leq \frac{(h+h^2)}{2}\|u''\|_{\infty}.
\end{equation}
\end{theorem}
\begin{prooff}
Here also, we begin by evaluating the $L^{1}$-norm of the derivative, that is $\|u'-u'_{I}\|_{0,1}$.
By the help of (\ref{T1})-(\ref{Remainder_TL}), setting $a=x_i, b=x_{i+1}$, we obtain that
$$
\int_{x_i}^{x_{i+1}}|u'(x)-u'_{I}(x)|\,dx = \int_{x_i}^{x_{i+1}}\bigg|u'(x) -\bigg( \frac{u'(x_i) + u'(x_{i+1})}{2n} + \frac{1}{n}\sum \limits_{k=1}^{n-1} u'\!\!\left(\!x_i + k\frac{h_i}{n}\right) + \epsilon_{n} \bigg) \bigg|\, dx.
$$
Writing now $u'(x)$ in the integral as
$$
\D u'(x) = \frac{1}{2n} u'(x) + \frac{1}{2n} u'(x) + \frac{1}{n}\sum_{k=1}^{n-1}u'(x),
$$
enables us to derive the following estimate:
\begin{equation}\label{I1I2I3}
\D \int_{x_i}^{x_{i+1}}|u'(x)-u'_{I}(x)|\,dx \leq I_1 + I_2 + I_3 + \frac{h_i^2}{8n}(M_2-m_2),
\end{equation}
where we set:
\begin{eqnarray}
I_1 & = & \frac{1}{2n}\int_{x_i}^{x_{i+1}}|u'(x)-u'(x_i)|\,dx,\nonumber \\[0.2cm]
I_2 & = & \frac{1}{2n}\int_{x_i}^{x_{i+1}}|u'(x)-u'(x_{i+1})|\,dx, \nonumber\\[0.2cm]
I_3 & = & \frac{1}{n}\sum\limits_{k=1}^{n-1} \int_{x_i}^{x_{i+1}}|u'(x)-u'(x'_k)|\,dx, \label{I3}
\end{eqnarray}
with: $\D x'_k= x_i +k\frac{h_i}{n}$. \sa
Implementing the same techniques we used in (\ref{IN1})-(\ref{IN6}), that is, based on Fubini's theorem, we get for $I_1$ and $I_2$:
\begin{eqnarray}
I_1 & \leq & \frac{1}{2n}\int_{x_i}^{x_{i+1}}(x_{i+1}-t)|u''(t)|\,dt \leq \frac{h_i^2}{4n}\,\|u''\|_{\infty}, \label{I1}\\[0.2cm]
I_2 & \leq & \frac{1}{2n}\int_{x_i}^{x_{i+1}}(t-x_i)|u''(t)|\,dt \leq \frac{h_i^2}{4n}\,\|u''\|_{\infty}. \label{I2}
\end{eqnarray}
Concerning the estimate of $I_3$, the treatment is based on the same principles but it has to be adapted. We proceed as follows. \sa
Having $x'_k$ which belongs to the open interval $]x_i,x_{i+1}[$, we first split the integral in (\ref{I3}) as
\begin{equation}\label{I3_1}
\D \int_{x_i}^{x_{i+1}}|u'(x)-u'(x'_k)|\,dx = \int_{x_i}^{x'_k}\bigg|\int_{x'_k}^{x} u''(t)\,dt\bigg|\,dx + \int_{x'_k}^{x_{i+1}}\bigg|\int_{x'_k}^{x} u''(t)\,dt\bigg|\,dx.
\end{equation}
Then, by using Fubini's theorem on each term separately, we get the two following inequalities:
$$
\D \int_{x_i}^{x'_k}\bigg|\int_{x'_k}^{x} u''(t)\,dt\bigg|\,dx \leq \int_{x_i}^{x'_k}(t-x_i)|u''(t)|\,dx,
$$
and
$$
\D \int_{x'_k}^{x_{i+1}}\bigg|\int_{x'_k}^{x} u''(t)\,dt\bigg|\,dx \leq \int_{x'_k}^{x_{i+1}}(x_{i+1}-t)|u''(t)|\,dx,
$$
which enables us to get  from (\ref{I3_1})
\begin{eqnarray}
\D \D \int_{x_i}^{x_{i+1}}|u'(x)-u'(x'_k)|\,dx & \leq & \int_{x_i}^{x'_k}\!(t-x_i)|u''(t)|\,dx + \!\int_{x'_k}^{x_{i+1}}\!(x_{i+1}-t)|u''(t)|\,dx, \nonumber \\[0.2cm]
& \leq & \bigg[\!\int_{x_i}^{x'_k}\!(t-x_i)\,dt + \!\int_{x'_k}^{x_{i+1}}\!(x_{i+1}-t)\,dt \bigg]\|u''\|_{\infty}, \nonumber \\[0.2cm]
\D & \leq & \!\!\!\! \D\bigg[\frac{(x'_k-x_i)^2}{2} + \frac{(x_{i+1}-x'_k)^2}{2} \bigg]\|u''\|_{\infty}. \label{I3_2}
\end{eqnarray}
Then, using that $x'_k\,\in\,]x_i,x_{i+1}[$, we can write it as a affine combination of $x_i$ and $x_{i+1}$, namely:
$$x'_k = tx_i +(1-t)x_{i+1}, (0<t<1).$$
Thus, (\ref{I3_2}) can be written as
\begin{equation}
\D \D \int_{x_i}^{x_{i+1}}|u'(x)-u'(x'_k)|\,dx \leq \ \frac{1}{2}\big(2t^2-2t+1\big)h_i^2\,\|u\|_{\infty} \leq \frac{h_i^2}{2}\,\|u''\|_{\infty},
\end{equation}
and by summing on $k$ between $0$ and $n-1$, (\ref{I3}) leads to:
\begin{equation}\label{I3_3}
\D I_3 \leq \bigg(\frac{n-1}{2n}\bigg)\,h_i^2\,\|u''\|_{\infty}.
\end{equation}
Finally, by the help of inequalities (\ref{I1}), (\ref{I2}) and (\ref{I3_3}), the estimate (\ref{I1I2I3}) gives
\begin{eqnarray}
\D \int_{x_i}^{x_{i+1}}|u'(x)-u'_{I}(x)|\,dx & \leq & \bigg[\frac{1}{4n} + \frac{1}{4n} + \frac{n-1}{2n} \bigg]\,h_i^2\,\|u''\|_{\infty} + \frac{h_i^2}{8n}(M_2-m_2),\nonumber \\[0.2cm]
\D & \leq & \frac{h_i^2}{2}\|u''\|_{\infty} + \frac{h_i^2}{8n}(M_2-m_2)\,.\label{I3_4}
\end{eqnarray}
Now, by summing on $i$ between $0$ to $N$, we get the following estimate for the $L^{1}$-norm of the derivative:
\begin{equation}\label{F1}
\D \|u'-u'_{I}\|_{0,1} \leq \frac{h}{2}\|u''\|_{\infty} + \frac{h}{8n}(M_2-m_2).
\end{equation}
Let us now evaluate the $L^1$-norm, that is $\|u-u_{I}\|_{0,1}$. \sa
Like in the previous proof, (see (\ref{Ineq0})-(\ref{Ineq00})), we can write, for all $x \in [x_i,x_{i+1}]$,
$$
\D |u(x)-u_{I}(x)| \leq \int_{x_i}^{x}|u'(t)-u'_I(t)|\,dt \leq \int_{x_i}^{x_{i+1}}|u'(t)-u'_I(t)|\,dt.
$$
Then, due to estimate (\ref{I3_4}), this inequality becomes
$$
\D |u(x)-u_{I}(x)| \leq \frac{h_i^2}{2}\|u''\|_{\infty} + \frac{h_i^2}{8n}(M_2-m_2).
$$
Now, integrating on $x$ which belongs to $[x_i,x_{i+1}]$ leads to
$$
\D \int_{x_i}^{x_i+1}|u(x)-u_{I}(x)|\,dx \leq \frac{h_i^3}{2}\|u''\|_{\infty} + \frac{h_i^3}{8n}(M_2-m_2).
$$
Finally, we sum over all values of $i$ between $0$ to $N$ and we get that
\begin{equation}\label{F2}
\D \|u-u_{I}\|_{0,1} \leq \frac{h^2}{2}\|u''\|_{\infty} + \frac{h^2}{8n}(M_2-m_2).
\end{equation}
In these conditions, the $W^{1,1}$-norm of the $P_1$-interpolation error can be obtained by adding inequalities (\ref{F1}) and (\ref{F2}), that is
\begin{equation}\label{F3}
\D \|u-u_{I}\|_{1,1} \leq \frac{h+h^2}{2}\|u''\|_{\infty} + \frac{h+h^2}{8n}(M_2-m_2), \forall n \in \N.
\end{equation}
By letting $n$ going to $+\infty$ in (\ref{F3}), we obtain estimate (\ref{Error_Estim_TL_Asymp}).
\end{prooff}
\begin{remark}
We can observe the impact Theorem \ref{Thmfirst} in various contexts, including  Lagrange finite element error estimate. Indeed,  let $u \in V$ denote the solution to a second order elliptic variational formulation \cite{RaTho82}, and $u_{h}$ its finite element approximation, where $V$ represents the solution space. According to C\'ea's lemma (see for instance \cite{Ern_Guermond}), to estimate the approximation error $\left\|u - u_{h} \right\|_{V}$, we need to select an element $v^{*}_{h}$ for which we can compute an estimate of $\left\|u - v^{*}_{h} \right\|_{V}$.\sa
A well-known and convenient choice consists in choosing $v^{*}_{h}$ as an interpolation polynomial of a specified degree. In this regard, the use of the new Taylor-like formula allows us to enhance the accuracy of the finite element solution for a given mesh size.
\end{remark}
\section{Enhancing finite differences schemes discretization}\label{Heat}\label{Section-newsecond}
\noindent In this section, we give another application of the new Taylor-like formula (\ref{T2}). We  first consider an implicit finite differences scheme to approximate 
the heat equation. Then, we will show how we can improve the error bound of the error estimate, without deteriorating neither the consistency error and the order of convergence, 
nor the stability.\sa
More precisely, we will derive with (\ref{T2}) a new implicit finite differences scheme, that is  first-order in time and second-order in space. In addition, it will be unconditionally stable, as the standard implicit finite differences scheme, but the upper bound of the error estimate will be two times smaller than the  standard one. \sa
To begin with, let us introduce the one dimensional heat equation defined as follows: consider a function $u(x,t)$ defined on $I\times[0,T]$, solution to:
\begin{equation}\label{Heat1}
\D \frac{\partial u}{\partial t}(x,t) = \frac{\partial^2 u}{\partial x^2}(x,t), \quad (x,t)\in I\times [0,T],
\end{equation}
where $I$ is a given open subset of $\R$, (bounded or not), and $T$ a given positive real number. \sa
We also introduce a constant time step $k\equiv\Delta t$ and the corresponding discrete time $t^{(n)}$, defined by: $t^{(n)}=nk, (n\in\N)$, as well as a constant mesh size $h\equiv\Delta x$ and the nodes of the mesh $x_j=jh, (j\in\Z)$. Finally, we denote by $\tilde{u}_j^{(n)}$ an approximation of solution $u$ to equation
(\ref{Heat1}) at a given point $(x_j, t^{(n)})$, that is, $\tilde{u}_j^{(n)} \simeq u(x_j, t^{(n)})$.\sa
Let us now remind the standard implicit scheme \cite{Euvrard}, that we denote \textbf{(FD$_1$)}:
\begin{equation}
\mbox{\textbf{(FD$_1$)}} \hspace{0.5cm} \frac{\tilde{u}_j^{(n+1)}-\tilde{u}_j^{(n)}}{k} = \frac{\tilde{u}_{j-1}^{(n+1)}-2\tilde{u}_j^{(n+1)}+\tilde{u}_{j+1}^{(n+1)}}{h^2}, \quad\forall (j,n)\in J\times\N,
\end{equation}
where $J$ denotes the total number of space nodes used in the space domain $I$. \sa
It is well known \cite{Euvrard} that the finite differences scheme \textbf{(FD$_1$)} is first-order in time and second-order in space. Indeed, introducing $m_2, M_2 \in \R$ 
the lower and upper bounds of the second derivative in time of $u$,  that is
$$
\D \forall t \in [0,T], \, m_{2} \leqslant \frac{\partial^2 u}{\partial t^2}(.,t) \leqslant M_{2}\,,
$$
and using the classical Taylor formula, we readily get that
$$
\D \frac{\partial u}{\partial t}(x,t^{(n)}) = \frac{u(x,t^{(n+1)})-u(x,t^{(n)})}{k} + \epsilon^{(T)}_{n},
$$
where
$$
|\epsilon^{(T)}_{n}| \leq  \frac{k}{2}\max(|m_2|,|M_2|).
$$
Finally, introducing $u_j^{(n)}\equiv u(x_j,t^{(n)})$,  the exact values of $u$ at the point $(x_j,t^{(n)})$, the well-known approximation of $\D \frac{\partial u}{\partial t}$ at $(x_j,t^{(n)})$ is given by
\begin{equation}\label{Approx_T}
\D \frac{\partial u}{\partial t}(x_j,t^{(n)}) \simeq \frac{u_j^{(n+1)}-u_j^{(n)}}{k}\,.
\end{equation}
Now, let us introduce another approximation of the first-order partial derivative $\D \frac{\partial u}{\partial t}$, for a solution $u$ to the equation (\ref{Heat1}) which is $C^2$ in time, by using the new Taylor-like formula (\ref{T2}).
\sa
So, let us set in (\ref{T2}) $a=t^{(n)}$ and $b=t^{(n+1)}$.  Then, for any given $x \in I$, we get: 
$$
\D u(x,t^{(n+1)}) = u(x,t^{(n)}) + \frac{k}{2}\bigg(\frac{\partial u}{\partial t}(x,t^{(n)}) + \frac{\partial u}{\partial t}(x,t^{(n+1)})\bigg) + k \epsilon_{n},
$$
where the remainder $\epsilon_{n}$ satisfies here
$$
|\epsilon_{n}| \leq  \frac{k}{8}(M_2-m_2).
$$
Then, by the help of (\ref{Heat1}), we obtain that
$$
\D \frac{\partial u}{\partial t}(x,t^{(n)}) = \frac{2}{k}\big(u(x,t^{(n+1)}) - u(x,t^{(n)})\big) - \frac{\partial^2 u}{\partial x^2}(x,t^{(n+1)}) - 2\epsilon_n,
$$
that leads to the following approximation of $\D \frac{\partial u}{\partial t}$ at the point $(x_j,t^{(n)})$:
\begin{equation}\label{Approx_TL}
\D \frac{\partial u}{\partial t}(x_j,t^{(n)}) \simeq \frac{2}{k}\big(u_j^{(n+1)} - u_j^{(n)}\big) - \frac{u_{j-1}^{(n+1)}-2u_j^{(n+1)}+u_{j+1}^{(n+1)}}{h^2}\,.
\end{equation}
This approximation is clearly of the first order in time and of the second order in space.\sa
Now, we want to compare the error bounds obtained by the approximations (\ref{Approx_T}) and (\ref{Approx_TL}) of the time partial derivative.  First, remark that the quantity involved in approximation (\ref{Approx_T}) is $\epsilon^{(T)}_{n}$, whereas $2\epsilon_n$ is involved in approximation (\ref{Approx_TL}). \sa
Moreover, assuming that $0\leq m_2 \leq M_2$ and setting $\Lambda \equiv M_2-m_2, (\Lambda \geq 0),$ we have:
$$
\D|2\epsilon_n| \leq \frac{k}{4}\Lambda\leq\frac{k}{2}\max\big(|m_2|,|m_2+\Lambda|\big)=\frac{k}{2}(m_2+\Lambda)\,.
$$
Consequently, approximation (\ref{Approx_TL}) leads to an upper bound necessary smaller than that of  (\ref{Approx_T}). \sa
The minimum gain that can be observed corresponds to $m_2=0$, and in this case, the error bound for approximation (\ref{Approx_TL}) is two times smaller than that of (\ref{Approx_T}). Moreover, as $m_2$ increases, the error bound related to (\ref{Approx_TL}) decreases at least twice.\sa
We are now in position to derive the main result of this section. So, for all $(j,n)\in J\times\N$, let us denote by \textbf{(FD$_2$)} the finite differences scheme defined by
\begin{equation}
\mbox{\textbf{(FD$_2$)}} \hspace{0.3cm}
-\frac{k}{2h^2}\,\tilde{u}_{j-1}^{(n+1)} + \bigg(1+\frac{k}{h^2}\bigg)\tilde{u}_{j}^{(n+1)} -\frac{k}{2h^2}\,\tilde{u}_{j+1}^{(n+1)} = \frac{k}{2h^2}\bigg(\tilde{u}_{j-1}^{(n)}+\tilde{u}_{j+1}^{(n)}\bigg) + \bigg(1-\frac{k}{h^2}\bigg)\tilde{u}_{j}^{(n)}.
\end{equation}
Hence, we have the following result:
\begin{theorem}
The implicit scheme \textbf{(FD$_2$)} is consistent with heat equation (\ref{Heat1}). It is first-order accurate in $k$ (time) and second-order accurate in $h$ (space), that is in $O(k)$ and $O(h)^{2}$. Moreover, for $m_2\geq 0$, the error estimate bound in $k$ is at least two times smaller that of the classical implicit finite differences scheme. Finally, when the open set of integration $I$ is equal to $\R$, the scheme \textbf{(FD$_2$)} is unconditionally stable.
\end{theorem}
\begin{prooff}
Consistency and order of convergence in space and in time were proved above when we introduced the scheme \textbf{(FD$_2$)}. \sa
To learn the stability of the scheme, and check the necessary stability condition, we introduce as usual the amplitude coefficient $A$ defined by $\tilde{u}_j^{(n)} = A^n e^{ipx_j}$, where $p$ is a given real parameter.
\sa
Then, replacing this expression of $\tilde{u}_j^{(n)}$ in the scheme \textbf{(FD$_2$)}, we get, after some elementary trigonometrical transformations, that the coefficient $A$ is equal to
\begin{equation}
\D A = \frac{1-X}{1+X}, \hs \mbox{ where } \hs X = \frac{2k}{h^2}\sin^2\big(\frac{ph}{2}\big).
\end{equation}
Given that $X\geq 0$, we obtain that $|A|\leq 1$, and the necessary stability condition is always fulfilled.\sa
Finally, when the open set of integration $I$ is equal to $\R$, the scheme \textbf{(FD$_2$)} is unconditionally stable, since the necessary condition becomes a sufficient one \cite{Euvrard}.
\end{prooff}

\begin{remark}
The application proposed here to the approximation of the heat equation has to be viewed as an example. Clearly, the new Taylor-like formula can be applied in various contexts and for different equations, depending on the specific requirements. Depending on their needs, readers can choose how to adapt formula (\ref{T1}) in a same manner as proposed here.
\end{remark}

\section{Conclusion}\label{D}
\noindent In this paper we have proposed a way to enhance interpolation and approximation error estimates. Based on a new Taylor-like formula, we show how it can be applied to improve first, the $P_1$-interpolation error estimate, then, to propose new numerical schemes yielding to smaller upper bounds in the corresponding error estimate.\\

\noindent In the first part, considering that numerous partial differential equations are well-posed within the general class of Sobolev spaces $W^{m,p}(\Omega), (m,p)\in\N^{*2}$, and particularly in the space $W^{1,1}(\Omega)$, we initially consider this framework, where we proved new interpolation error estimates.\\

\noindent  In a second part, we have considered another application of the new Taylor-like formula. Turning out our attention to the analysis of a classical implicit finite differences scheme used to discretize the heat equation, we showed how we can enhance the upper bound of the error estimate, while maintaining the consistency error and the order of convergence at the same level, together with the stability as well.\\

\noindent This new first-order Taylor-like formula could also find applications in various other contexts. For example, we can develop new schemes for different types of partial differential equations or ordinary differential equations, as well as tackle more general problems within the Sobolev spaces $W^{m,p}(\Omega)$ when
$p>1$.\\

\noindent \textbf{\underline{Homages}:} The authors want to warmly dedicate this research to pay homage to the memory of Professors Andr\'e Avez and G\'erard Tronel who largely promote the passion of research and teaching in mathematics of their students.
\end{document}